\documentclass[a4paper,12pt, legno]{article}
\usepackage{amsmath}
\usepackage{amsfonts}
\usepackage{amsthm}
\usepackage{mathrsfs}
\usepackage[T1]{fontenc}
\usepackage{amsmath}
\usepackage{amsfonts}
\usepackage{amsthm}
\usepackage{mathrsfs}
\usepackage{amssymb}
\usepackage{pst-node}
\usepackage{tikz-cd}

\newtheorem{thm}{Theorem}

\newtheorem{lem}[]{Lemma}

\newtheorem{cor}{Corollary}

\newtheorem{rem}[]{Remark}

\newcommand{\Rset}{\mathbb{R}}
\newcommand{\Cset}{\mathbb{C}}

\begin{document}

\title{Remarks on  the tangent function  from an analytic and a probability  point of view.}
\author{Zbigniew J.  Jurek  (University of Wroc\l aw}
\date{February 15, 2021}
\maketitle
\maketitle

\begin{quote}
\textbf{Abstract.}
In this note we prove that the function $\tan(1/it)$  is a Voiculescu transform of a free-infinitely divisible distribution, that is, it admits the integral representation associated with  Pick functions. Moreover, we found its "counterpart" in the classical infinitely divisible measures expressed by a series of Rademacher variables.

\medskip



\medskip
\emph{Abbreviated title: On a tangent function from probability point of view.}

\end{quote}

\medskip
 \underline{Author's address}:

Institute of Mathematics, University of Wroc\l aw, Pl. Grunwaldzki 2/4,

50-384 Wroc\l aw, Poland;

Email:  zjjurek@math.uni.wroc.pl ;   www.math.uni.wroc.pl/$\sim$zjjurek

ORCID Number: 0000-0002-8557-2854

\medskip
\medskip
\medskip
In this note, by quite elementary arguments, we prove  that  the tangent  function  $\tan(1/it), t\neq 0$, is  a transform of a free-infinitely divisible probability measure. Thus it admits  the integral representation (Pick function restricted to the imaginary line); for instance cf. Bondensson (1992), p. 20.

However, if we know that both $1/z$ and $\tan(z)$ are Pick functions, and the property of being Pick function is closed under  a superposition, we immediately  conclude that $\tan(1/z)$ is a Pick function as well and hence we infer the representation in question; see  F. Gesztesy and E. Tsekanovskii (2000) and references therein. For some explicit examples, see there in  Section A, p.121.

Our first aim here is not to use this structural property of Pick functions and the second one is to find for the tangent function its   $"$counterpart part$"$  in the  classical infinite divisibility using  an isomorphism given in Jurek (2016). More explicitly,  we do it via a random integral mapping $\mathcal{K}$ from Jurek (2007) and series representation of a hyperbolic tangent function; cf. Theorem 1.

In Jurek (2020) an analogy between different notions of infinite divisibility was studied and there, in Section 4, are given  explicit examples of Pick functions -- free-infinite divisible transforms -- related to the hyperbolic and Laplace (double exponential)  characteristic functions.

\medskip
\medskip
\textbf{1. Notations and  terminology.}

For a finite Borel measure $m$ on the real line $\Rset$ and a Borel function $g$ on a positive halfline $[0,\infty)$ we define \emph{a characteristic function} $\phi_m(t)$ and \emph{ a Laplace transform} $\mathfrak{L}[g;w]$  as follows
 \begin{equation}
\phi_m(t):= \int_{\Rset} e^{itx} m(dx), \  t \in \Rset; \ \ \ \mathfrak{L}[g(x); w]:=\int_{0}^\infty g(x)\,e^{-w x}dx, \  w>0.
\end{equation}
A probability measure  $e(m):= e^{-m(\Rset)}\sum_{k=0}^\infty\frac{m^{\ast k}}{k!}$, where $\ast$ denotes a convolution of measures, is called \emph{a compund Poisson measure}. We will write $\mu=[a,\sigma^2, M]$, if  a probability measure $\mu$ is \emph{infinitely divisible} (that is: for each $k \ge 2$ there exist measure $\mu_k$ such that $\mu_k^{\ast k}=\mu$)  and
\begin{equation}
\phi_{\mu}(t)=  \exp[ita-\frac{1}{2} \sigma^2 t^2+\int_{\Rset\setminus{\{0\}}}(e^{itx}-1-\frac{itx}{1+x^2})M(dx)], \  t \in \Rset;
\end{equation}
where  $a\in \Rset$ (a shift),  $\sigma^2$  (a variance of a Gaussian part) and $M$ is a  measure  satisfying a condition  $\int_{\Rset} \min (1, x^2)M(dx)<\infty$ ( so called \emph{a L\'evy measure}). The triplet $a,\sigma^2, M$ is uniquely determined by $\mu$
and formula (2)  is  called \emph{the L\'evy-Khintchine representation}. By $(ID, \ast)$ we denote a convolution semigroup set of all infinitely divisible probability measures.
Note that compound Poisson measures $e(m)$ are  infinitely divisible and their characteristic functions are equal to  $\phi_{e(m)}(t)=\exp\int_{\Rset}(e^{itx}-1)m(dx)$.

A stochastic process $Y(t), t\ge 0,$ on a probability space $(\Omega, \mathcal{F},P)$, starting from zero ($Y(0)=0$ with P.1), with stationary and stochastic  independent increments, and continuous in probability is called \emph{a L\'evy process}. Hence  all increments $Y(t)-Y(s), t>s\ge 0$ have ID distributions and a probability distribution $\mu_t$ of random variable $Y(t)$ is equal to $\mu_1^{\ast t}$. And conversely, each  $\mu\in ID$ can be inserted into a L\'evy process $Y_{\mu}(t),  t\ge 0$, such that probability distribution of  $Y_{\mu}(1)=\mu$.

Finally,  below we need an random integral mapping
$$ (ID,\ast)\ni \mu \to \mathcal{K}(\mu)=  \mathcal{L}(\int_0^\infty s dY_{\mu}(1-e^{-s})),$$
where $\mathcal{L}(V)$  denotes a probability distribution of a random variable $V$, that was introduced in Jurek (2007).

Our investigations of the tangent function  have the origin in the following series representations:
\begin{equation}
\tanh s =  2  s \sum_{k=1}^\infty \frac{1}{((2k-1)\pi/2)^2 +s^2}= s\int_{\Rset} \frac{1}{x^2 +s^2}\,(\sum_{n \in Z,\, odd}\delta_{n\pi/2 } (dx))
\end{equation}
where $\delta_a$ is a Dirac measure (point-mass measure concentrated at $a$); cf. Gradshteyn and Ryzhik (1994), formula \textbf{1.421}(2).

\medskip
\textbf{2. Results.}

\medskip
Here we prove the following analytic and probabilistic properties of  a tangent function $\tan(1/it)$ :
\begin{thm}
(a) \ \ Let's  for  $k\in N :=\{1,2 ,3,...\}$  define real constants
 $$c_k:=( (2k-1)\pi/2)^{-1},  c_{-k}:=(2(-k)+1)\pi/2)^{-1}$$
 and   define discrete   measures
\begin{equation}
 m(dx):= \sum_{k \in Z}\frac{x^2}{1+x^2}\delta_{c_k}(dx);  \  \ \   \   M(dx):=\frac{1+x^2}{x^2}m(dx);
\end{equation}
where $\delta_c$ denotes the Dirac measure. Then $m$ is a finite measure with a mass  $m(\Rset)=\tanh(1) \approx 0.76159$, and
\begin{equation}
\   \int_{\Rset}\frac{1+itx}{it-x}m(dx)=  \tan(1/it) \ \ \ (\equiv \ -i \tanh(1/t)), \ t\neq 0;
\end{equation}
(b) \  \ Since $\int_{(\Rset\setminus{{0}})} (x^2 \land1)M(dx)<\infty$ therefore   $\mu=[0,0,M]$ represents a classical infinitely divisible
measure  with zero shift,  zero Gaussian part and L\'evy (spectral) measure $M$.  Let $\phi_{\mu}$ be its characteristic function.
Then
\begin{equation}
 it^2\, \mathfrak{L}[\log \bar{\phi_{\mu}}(s); t]  = it \log(\phi_{\mathcal{K}(\mu)}(- t^{-1}))=
 \tan(1/it) \ (\equiv - i \tanh(1/t)),
\end{equation}
where $\bar{\phi_{\mu}}$  is a complex conjugation  and $\mathfrak{L}[g; t]$ denotes a Laplace  transform of a function $g$ at $t$.

(c) for a measure $m$ and its characteristic function $\phi_m$ we have
\begin{equation}
(w^2-1)\mathfrak{L}[\phi_m(x)-\tanh(1)\cosh(x);w]= - i\tan(\frac{1}{iw}), \ \ (\equiv - \tanh(1/w))
\end{equation}
where $\mathfrak{L}$ stands for  a Laplace transform.

(d) \ \ A measure $\mu=[0,0,M]\in (ID,\ast)$,  a "counterpart" of Pick function $\tan(1/it)$  (free-infinite divisible transform), is a probability distribution of a variable  $X$  with the following series representation
\begin{equation}
X=\sum_n c_n Y_n<\infty \ \mbox{converges a.s., in distribution and in $L^2$ },
\end{equation}
where $Y_1,Y_2,...$ are i.i.d. copies of  a random sum $Y=\sum_{j=0}^{N_2}r_j$ with $N_2, r_1,r_2,...$ independent, where
$N_2$ has Poisson distribution  with the parameter $2$ and $r_j$  are Rademacher variables with $P(r_j = \pm 1)=1/2$.
\end{thm}

Since a function $\tan(z)$ is analytic one, from (5) via an extension arguments, we have that
\begin{cor}
For  $z \in \Cset\setminus{\{0\}}$
\[
\tan(1/z)=  \int_{\Rset}\frac{1+zx}{z-x}m(dx), \ \mbox{where} \ \ m(dx):=\sum_{n \in Z, odd}
\frac{x^2}{1+x^2}\delta_{ (n\pi/2)^{-1}} (dx).
\]
\end{cor}
This representation is also in Gesztesy and Tsekanovskii (2000), Example A.10 and A.11 on p.122 but  obtained by different complex analysis arguments.  Also note that  their examples  A.8 and A.9  correspond to stable probability measures with an exponent r. Also $\tan(1/z)$,  and much more of free-probability, can be found   in  Ejsmont and Lehner (2020).

\medskip
The following  lemma gives  a probability distribution of  the variable $Y$  that appears in the part (d) of Theorem 1  and may be  of some independent interest as well.
\begin{lem}
Let $N_2, r_1,r_2,..., r_k,...$ be  independent variable where $N_2$ has a Poisson distribution with a parameter 2 and $r_j=\pm 1$ with probability 1/2 are independent  Rademacher variables. Then
$$
Y:=\sum_{j=1}^{N_{2}} \, r_j \stackrel{d}{=}S_+ - S_{-}, \ \mbox{where $S_{\pm}$ are indpendent Poisson with parameter 1}.
$$
Moreover, for $k=0,1,2,...$ we have
\[
P(Y=k)=P(Y=-k)=e^{-2}\sum_{j=0}^\infty \frac{1}{j!} \frac{1}{(k+j)!}\equiv e^{-2} I_k(2);  \mbox{ (a Bessel function)}.
\]
In particular, $ I_0(2)+2 \sum_{k=1}^\infty I_k(2)=e^2$.
\end{lem}
(cf. also Remark 2, below).
\begin{rem}
\emph{ The above Lemma 1,  is a particular case of so called Skellam distributions  that are defined as a law of a difference  $T_1-T_2$ of  independent Poisson variables with a parameters $\lambda_1, \lambda_2$, respectively.}
\end{rem}

\medskip
\textbf{3. PROOFS.}

\medskip
\emph{Proof of Theorem 1.}

\underline{Part (a).}  Noticing that for $k\in N$  we  have $c_{-k}=-c_k$  and using (4) we
\begin{multline*}
\int_{\Rset}\frac{1+itx}{it-x}m(dx)= \sum_{k \in Z\setminus{\{0\}}} \frac{1+itc_k}{it - c_k}\frac{ c_k^2}{1+ c_k^2}= \sum_{k \in N} [ \frac{1+itc_k}{it - c_k}+  \frac{1+itc_{-k}}{it - c_{-k} }]\frac{c_k^2}{1+c_k^2}\\
= -2 it \sum_{k \in N}\frac{1+c_k^2}{t^2+c_k^2}\frac{c_k^2}{1+c_k^2}= - 2it \sum_{k\in N} \frac{1}{t^2}\frac{1}{c_k^{-2} +1/t^2} \  \ \  (\mbox{by(3)})\\=
 - i\, 2\,\frac{1}{t} \sum_{k\in N}\frac{1}{((2k-1)\pi/2)^2 + 1/t^2}= - i \tanh(1/t)= \tan(1/it),
\end{multline*}
which completes a proof of (4). Furthermore, using   the one to the last equality above  for $t=1$, we conclude
$$ m(\Rset)= \sum_{k \in Z}\frac{c_k^2}{1+c_k^2} = 2\, \sum_{k\in N}\frac{1}{((2k-1)\pi/2)^2 + 1}=\tanh(1),$$
and this completes a proof of a part $ (a)$.

\medskip
\underline{Part (b).} In Jurek (2007) a definition of the mapping $\mathcal{K}$ is given  by a formula (12) and using Corollary 6 we get first  equality in (6).  (Also compare  Jurek (2016)). Consequently, if
$$ (ID,\ast)\ni \mu \to \mathcal{K}(\mu)=  \mathcal{L}(\int_0^\infty s dY_{\mu}(1-e^{-s}))$$
then its  characteristic function is of   the  form:
$$  \log (\phi_{\mathcal{K}(\mu)}) (t)   =\int_0^\infty\log(\phi_{\mu}(st))e^{-s}ds, \  \ ( \mbox{by Corollary 4,  in Jurek(2007)}).$$
Specifying above for $\mu=[0,0,M]$,  with a symmetric $M$, we  get
\begin{multline*}
 \log (\phi_{\mathcal{K}([0,0,M])}) (t) =\int_0^\infty\int_{\Rset\setminus{\{0\}}}(\cos (tsx)-1))M(dx)e^{-s}ds\\  =
\int_{\Rset\setminus{\{0\}}}[\int_0^\infty(\cos(stx)-1)e^{-s}ds]M(dx) =\int_{\Rset\setminus{\{0\}}}[\frac{1}{t^2x^2 +1}-1]M(dx) \\ = - t^2\int_{\Rset\setminus{\{0\}}}\frac{x^2}{t^2x^2 +1}M(dx)  =
-t^2 \sum_{k\in Z}\frac{c_k^2}{t^2c_k^2  +1} =-t^2 \sum_{k\in Z}\frac{1}{t^2 + c_k^{-2}} = \\  - t  \tanh(t)= it \tan(it). \qquad \qquad  \qquad  \qquad  \qquad  \quad \quad
\end{multline*}
Substituting above $t$ by $-1/t$ we get formula (6).

\medskip
\underline{Part (c).} Using a measure $m$, a formula (3),  part (a) of Theorem 1 and Laplace transforms of a  cosine and hyperbolic cosine functions  we get
\begin{multline*}
\mathfrak{L}[\phi_m(x)-\tanh(1) \cosh(x); w]\\ = \int_0^\infty[\sum_{n\in Z, odd}\cos(w c_n)\frac{c_n^2}{1+c_n^2} -\cosh(x)\sum_{n\in Z, odd}\frac{c_n^2}{1+c_n^2}\,]\,e^{-wx}dx \\=
\sum_{n \in Z, odd}[ \frac{w}{w^2+c_n^2}\frac{c_n^2}{1+c_n^2} - \frac{w}{w^2-1}\frac{c_n^2}{c_n^2 +1}]\\ =
 \frac{w}{w^2-1} \sum_{n \in Z, odd} (\frac{w^2-1}{w^2+c_n^2} -1)
\frac{c_n^2}{1+c_n^2}= -  \frac{1}{w^2-1}\, \frac{1}{w} \sum_{n \in Z, odd} \frac{c_n^2}{1+c_n^2w^{-2}} \\ = (w^2-1)^{-1} (-1) \tanh (\frac{1}{w}) = (w^2-1)^{-1}\, \,(-i)\, \tan(\frac{1}{ iw}),
\end{multline*}
which completes a proof of a part (c).

\medskip
\underline{Part(d).}
First, note that $\cos t$ is a characteristic function of Rademacher variables $r_j$. Second, variable $Y$ has a compound Poisson distribution with a characteristic function $\exp[2(\cos t-1)]$.Third, since $\mathbb{E}[Y]=0, Var[Y]= 2$ therefore $\mathbb{E}[X]=0$ and $Var(X)=\sum_n Var(c_nY_n)= 2  \sum_n c_n^2<\infty$. All in all, by Kolmogorov's Three Series Theorem we conclude that  series (8) of independent variables converge in all three modes.

Finally, we have
\[
\phi_X(t)=\prod_n \exp[2(\cos(c_nt)-1)]= \exp\int_{\Rset}(\cos (tx)-1)M(dx)=\phi_{\mu}(t),
\]
which proves (8) and this completes a proof of Theorem 1.

\medskip
\medskip
\emph{Proof of Lemma 1.}

Let
$S_+$ and $S_{-}$ denote a number of $(+1)$   and $-1$ among $r_1,r_2,..., r_{N_2}$, respectively. Then
$$[S_+=k] \ \mbox{iff} \  [ (N_2\ge k) \ \mbox{and} \  ( \sum_{j=1}^{N_2} r_j = k (+1) + (N_2-k)(-1)= 2k-N_2].$$
Hence  using conditioning arguments ($S_+$ and $S_-$ are $N_2$ conditionally independent) we get
\begin{multline*}
P(S_+=k)=P( (  \sum_{j=1}^{N_2} r_j =2k-N_2)\cap (N_2\ge k))=\mathbb{E}[ 1_{\{N_2\ge k\}} \binom{N_2}{k}2^{- N_2} ]\\ =
\sum_{j=k}^\infty \binom{j}{k} \frac{1}{2^j}\ e^{-2} \frac{2^j}{j!} =
e^{-2} \sum_{j=k}^\infty \frac{j!}{k!\, (j-k)!} \frac{1}{j!}=\frac{1}{k!}e^{-2}\sum_{n=0}^\infty\frac{1}{n!}=\frac{1}{k!}e^{-1},
\end{multline*}
and therefore $S_+$ has a Poisson distribution with a parameter 1.

By a similar argument $S_{-}$ has also Poisson distribution withe a parameter 1 and is independent of $S_+$.  Finally,
\begin{multline*}
P(Y=k)=P(S_+ - S_{-}=k)= P(S_+- S_{-}=-k)     \\ = P((S_+=k+S_{-})\cap (\bigcup_{j=0}^\infty(S_{-}= j))
=\sum_{j=0}^\infty P((S_+=k+j)\cap(S_{-}=j))  \\ =\sum_{j=0}^\infty P(S_+=k+j)\,P(S_{-}=j)
=\sum_{j=0}^\infty e^{-1}\frac{1}{(k+j)!}e^{-1}\frac{1}{j!} \equiv e^{-2}I_{k}(2).
\end{multline*}
(See Remark 3 below  for a definition of  a modified Bessel function $I_k(z)$.)

\begin{rem}
\emph{
Since $I_k(2)$ are a probability distribution therefore
\begin{equation}
2\sum_{k=1}^\infty I_k(2)+ I_0(2)=e^2, \ \ \mbox{where} \ I_k(2):= \sum_{j=0}^\infty \frac{1}{j!} \frac{1}{(k+j)!}, \ \ k\ge 0.
\end{equation}
 But here is a straightforward calculation of that fact :
\begin{multline*}
e^2= \sum_{k, l =0}^\infty \frac{1}{k!} \frac{1}{l!}= \sum_{k=0}^\infty (\frac{1}{k!})^2  + 2 \sum_{k=0}^\infty\sum_{l=k+1}^\infty \frac{1}{k!} \frac{1}{l!}= I_0(2)+2 \sum_{k=0}^\infty\sum_{j=1}^\infty \frac{1}{k!} \frac{1}{(k+j)!}\\ =I_0(2)+2 \sum_{k=0}^\infty \big[\frac{1}{k!} ( \, \sum_{j=0}^\infty  \frac{1}{(k+j)!}- \frac{1}{k!}\,)\big]= I_0(2)+2 \sum_{j=0}^\infty  \sum_{k=0}^\infty \big( \frac{1}{k!}  \frac{1}{(k+j)!}- (\frac{1}{k!})^2\,\big)
\\= I_0(2) + 2 \sum_{j=0}^\infty ( I_j(2)-I_0(2))=  I_0(2) + 2 \sum_{j=1}^\infty  I_j(2);  \qquad \qquad  \qquad
\end{multline*}
which completes a proof of (9)}.
\end{rem}

\begin{rem}
\emph{
Recall that
 $$ I_{\nu}(z):= \sum_{j=0}^\infty\frac{1}{j! \Gamma(\nu +j+1)} (\frac{z}{2})^{\nu+2j},
 \  \  z\in\Cset,$$  is \emph{a modified  Bessel function of the first kind}; cf.  Gradshteyn and Ryzhik(1994),formula \textbf{8.445}.
Moreover, it has the following  integral representation
$$
I_{\nu}(z)=\frac{(\frac{z}{2})^\nu}{\Gamma(\nu +1/2)\Gamma(1/2)}\int_{-1}^1 (1-t^2)^{\nu-1/2}e^{\pm  z t}dt;  \ \Re(\nu+1/2)>0;
$$
Gradshteyn and Ryzhik (1994), formula \ \textbf{8.431(1)}.
\newline
Since $\Gamma(1/2)=\sqrt{\pi};  \ \Gamma(k+1/2)=\sqrt{\pi}2^{-k}(2k-1)!! $ then
$$
I_k(2)=\frac{2^k}{(2k-1)!! \, \pi}\int_{-1}^1 (1-x^2)^{k-1/2}e^{- 2 x}dx,  \ \ \mbox{for} \ \ k \in N,
$$
that might be more useful for an explicit calculations.}
\end{rem}

\medskip
\textbf{Acknowledgment.}
This note was prompted by a very useful discussion  with  Dr W. Ejsmont (University of Wroclaw) who shared with Author his work in progress  and   provided  a reference  to Gesztesy and Tsekanovski (2000).

\medskip
\medskip
REFERENCES.

\medskip
[1] \  L. Bondesson (1992), \emph{Generalized gamma convolutions and related}

\emph{ classes of distributions and densities},
Springer-Verlag, New York; Lecture

Notes in Statistics, vol. 76.

\medskip
[2]  \ W. Ejsmont and F. Lehner (2020),  The free tangent law,

arXiv:2004.02679.

\medskip
[3] \ F. Gesztesy and E. Tsekanovskii (2000),  On  matrix-valued Herglotz

functions, \emph{Math. Nachr.} vol. \textbf{218} , pp. 61-138.

\medskip
[4] \ I. S. Gradeshteyn and I. M. Ryzhik (1994), \emph{ Tables of integrals, series,}

\emph{and products}, Academic Press, San Diego, New York, Boston; $5^{th}$ Edition.

\medskip
[5]  \ Z. J. Jurek (2007),  Random integral representations for

free-infinitely divisible  and tempered stable distributions,

\emph{Stat.\& Probab. Letters},
    \textbf{77} no. 4 ,  pp. 417-425.

Also on:  arXiv:1206.3044,  [math\  PR], 14 June 2012.

\medskip
[6] \ Z. J. Jurek (2016),  On a method of introducing free-infinitely

divisible probability measures, \emph{Demonstratio Mathematica},  \textbf{49},

 no 2,   pp. 236-251.

Also on:   arXiv: 1412.4445, [math \ PR], 14  December 2014.

\medskip
[7]\ Z. J. Jurek (2020), On a relation between classical and free infinitely

divisible transforms, \emph{Probab. Math. Stat.}, vol. 40, Fasc. 2, pp. 349-367.

Also on:  arXiv:1707.02540 [math PR], 9 July 2017.

\end{document}